\documentstyle[12pt,amssymb]{amsart}

\renewcommand{\(}{\left(}
\renewcommand{\)}{\right)}
\renewcommand{\[}{\left(}
\renewcommand{\]}{\right)}
\newcommand{\<}{\langle}
\renewcommand{\>}{\rangle}
\newcommand{\p}{\rho}

\newcommand{\im}{{\mathrm Im}}
\renewcommand{\S}{{\mathcal S}}
\renewcommand{\SS}{S^\S}
\newcommand{\RA}{{\mathcal R}}

\newcommand{\tr}{{\mathrm tr}}
\newcommand{\R}{{\mathbb R}}

\newcommand{\norm}[1]{\left\lVert#1\right\rVert}
\renewcommand{\bar}{\overline}
\renewcommand{\phi}{\varphi}
\newcommand{\e}{\varepsilon}

\theoremstyle{plain}
\newtheorem{thm}{Theorem}[section]
\newtheorem{lem}[thm]{Lemma}

\newtheorem{cor}[thm]{Corollary}

\newtheorem*{MT*}{Main theorem}
\newtheorem*{MC*}{Main corollary}

\theoremstyle{definition}
\newtheorem{defn}[thm]{Definition}

\theoremstyle{remark}

\newtheorem*{PMT*}{Proof of the main theorem}

\title[Spectral synthesis]{A non-commutative Sobolev estimate\\
and its application to spectral synthesis}\thanks{
Based in part on the author's doctoral thesis (University of Chicago)
written under the direction of Professor Tim Steger.
}
\author{M. K. Vemuri}
\address{Chennai Mathematical Institute,
SIPCOT Information Technology Park,
Navalur Post,
SIRUSERI 603 103
India
}
\email{mkvemuri\@cmi.ac.in}



\begin{document}

\begin{abstract}
In [M. K. Vemuri, {\em Realizations of the canonical representation,}
submitted], it was shown that the spectral synthesis problem for
the Alpha transform is closely related to the problem of
classifying realizations of the canonical representation (of the
Heisenberg group).  In this paper, we show that discrete sets are sets
of spectral synthesis for the Alpha transform.
\end{abstract}

\maketitle

\tableofcontents\newpage

\section{Introduction}

For $p \in [1, \infty]$, let $S^p$ denote the Schatten class of 
$p^{\text{th}}$ power traceable operators on $L^2(\R)$.  For
$x,y\in\R$, let $T_x$ and $M_y$ denote the translation and modulation
operators on $L^2(\R)$, i.e. for $s\in L^2(\R)$
\begin{align*}
(T_x s)(t)=& s(t-x) \qquad\text{and} \\
(M_y s)(t)=& e^{-2\pi iyt} s(t).
\end{align*}
If $(x_1,y_1) \in \R^2$ and $X \in S^p$, set
$$(x_1,y_1) \cdot X = T_{y_1} M_{-x_1} X M_{x_1} T_{-y_1}.$$
If $q \in L^1(\R^2)$ and $X \in S^p$, set
$$q \cdot X = \iint q(x_1,y_1) (x_1, y_1) \cdot X \,dx_1 \,dy_1.$$
Then $S^p$ becomes an $L^1(\R^2)$-module.  Note that the previous integral
exists by Lemma \ref{L:LOI1}.

\begin{MT*}
If $X \in S^1$, $\tr(X)=0$ and $\e > 0$, then there exists
$\p \in L^1(\R^2)$ with $\hat\p=1$ on a neighborhood of $(0,0)$ such that
$$ \norm{\p \cdot X}_{S^1} < \e.$$
\end{MT*}

\begin{defn}
If $X \in S^1$, the Alpha transform of $X$ is the function on
$\R^2$ defined by
$$\alpha(X)(x,y)=\tr(T_x M_y X).$$
\end{defn}

The Alpha transform is related to the module structure on
$S^1$ in the same way that the classical Fourier transform is
related to convolution.  In fact, we have the following lemma.

\begin{lem}
If $q \in L^1(\R^2)$ and $X \in S^1$, then 
$\alpha(q \cdot X) = \hat q \alpha(X)$.
\end{lem}

\begin{pf}
Note first that
\begin{align*}
\alpha((x_1,y_1)\cdot X)
=& \tr(T_x M_y T_{y_1} M_{-x_1} X M_{x_1} T_{-y_1}) \\
=& e^{-2\pi i yy_1} \tr(T_{y_1} T_x M_{-x_1} M_y X M_{x_1} T_{-y_1}) \\
=& e^{-2\pi i (yy_1+xx_1)} \tr(T_{y_1} M_{-x_1} T_x M_y X M_{x_1} T_{-y_1}) \\
=& e^{-2\pi i (yy_1+xx_1)} \tr(T_x M_y X) \\
=& e^{-2\pi i (yy_1+xx_1)} \alpha(X)(x,y).
\end{align*}
The result now follows by integration.
\end{pf}

By standard methods, the main theorem leads to the following
corollary, which may be viewed as saying that discrete sets are sets
of spectral synthesis (more precisely, C-sets) for the Alpha
transform.

\begin{MC*}
Let $D \subseteq \R^2$ be a discrete set.  Let $\e > 0$.
If $X \in S^1$ and $\alpha(X)$ vanishes on $D$, then there exists
$Y \in S^1$ such that $\alpha(Y)$ vanishes on a neighborhood
of $D$ and $\norm{X-Y}_{S^1} < \e$.
\end{MC*}

\section{Some lemmas on integration}

In this section, we prove some slight extensions of \cite[Theorem 3.27]{Rudin-FA}.

\begin{defn}
Let $f$ be a continuous function defined on $\R^2$.  We say that $f$ is
{\em rapidly decreasing} and write $f \in \RA(\R^2)$ if for all
positive integers $n$, the function $(x,y) \mapsto (1+x^2+y^2)^n f(x,y)$ is
bounded.  The best constants in the bounds give a countable family of
norms which turn $\RA(\R^2)$ into a Frechet space.
\end{defn}

Let ${\mathcal T}$ be a Frechet space defined by a countable family
of {\em norms}.

\begin{defn}
A function $\Phi:\R^2 \to {\mathcal T}$ is
{\em bounded} if for each norm $\sigma$, the function $\sigma \circ \Phi$ is
bounded.
\end{defn}

\begin{defn}
A function $\Phi:\R^2 \to {\mathcal T}$ is
{\em polynomially bounded} if for each norm $\sigma$, the function $\sigma \circ \Phi$ is
dominated by a polynomial.
\end{defn}

\begin{lem}\label{L:LOI1}
Let $\Phi:\R^2 \to {\mathcal T}$ be a bounded continuous function.  Let
$f \in L^1(\R^2)$.  Then
$$\iint f(x,y) \Phi(x,y) \,dx\,dy $$
exists.
\end{lem}

\begin{pf}
There is a sequence $\{f_k\}$ of compactly supported continuous functions such that $f_k \to f$ in $L^1(\R^2)$.
For each $k$, the integral
$$ I_k = \iint f_k(x,y) \Phi(x,y) \,dx\,dy $$
exists by \cite[Theorem 3.27]{Rudin-FA}.
For each norm $\sigma$,
\begin{align*}
\sigma(I_k - I_j) \le & \iint |f_k(x,y)-f_j(x,y)| \sigma(\Phi(x,y)) \,dx\,dy \\
                  \le & \sup \sigma(\Phi(x,y)) \norm{f_k - f_j}_1 \\
                  \to & 0
\end{align*}
as $k,j \to \infty$.  So $\{I_k\}$ is a Cauchy sequence.  Since ${\mathcal T}$
is Frechet, there exists $I \in {\mathcal T}$ such that
$$ \lim_{k \to \infty} I_k = I.$$
Let $\Lambda$ be a continuous linear functional on ${\mathcal T}$.  Then
\begin{align*}
\Lambda(I) = & \lim_{k \to \infty} \Lambda(I_k) \\
           = & \lim_{k \to \infty} \Lambda(\iint f_k(x,y) \Phi(x,y) \,dx\,dy) \\
           = & \lim_{k \to \infty} \iint f_k(x,y) \Lambda(\Phi(x,y)) \,dx\,dy \\
           = & \iint f(x,y) \Lambda(\Phi(x,y)) \,dx\,dy \qquad\text{(by the dominated convergence theorem)} \\
           = & \iint \Lambda(f(x,y) \Phi(x,y)) \,dx\,dy.
\end{align*}
Therefore, $\iint f(x,y) \Phi(x,y) \,dx\,dy = I$.
\end{pf}

\begin{lem}\label{L:LOI2}
Let $\Phi:\R^2 \to {\mathcal T}$ be a polynomially bounded continuous function.  Let
$f \in \RA(\R^2)$.  Then
$$\iint f(x,y) \Phi(x,y) \,dx\,dy $$
exists.
\end{lem}

\begin{pf}
There is a sequence $\{f_k\}$ of compactly supported continuous functions such that $f_k \to f$ in $\RA(\R^2)$.
For each $k$, the integral
$$ I_k = \iint f_k(x,y) \Phi(x,y) \,dx\,dy $$
exists by \cite[Theorem 3.27]{Rudin-FA}.
For each norm $\sigma$, there exist $C$ and $n$ such that
$$ \sigma(\Phi(x,y)) \le C(1+x^2+y^2)^n. $$
Thus
\begin{align*}
\sigma(I_k - I_j) \le & \iint |f_k(x,y)-f_j(x,y)| \sigma(\Phi(x,y)) \,dx\,dy \\
                  \le & C \iint (1+x^2+y^2)^n |f_k(x,y) - f_j(x,y)| \,dx\,dy \\
                  \to & 0
\end{align*}
as $k,j \to \infty$.  So $\{I_k\}$ is a Cauchy sequence.  Since ${\mathcal T}$
is Frechet, there exists $I \in {\mathcal T}$ such that
$$ \lim_{k \to \infty} I_k = I.$$

Let $\Lambda$ be a continuous linear functional on ${\mathcal T}$.  Then
\begin{align*}
\Lambda(I) = & \lim_{k \to \infty} \Lambda(I_k) \\
           = & \lim_{k \to \infty} \Lambda\(\iint f_k(x,y) \Phi(x,y) \,dx\,dy\) \\
           = & \lim_{k \to \infty} \iint f_k(x,y) \Lambda(\Phi(x,y)) \,dx\,dy \\
           = & \iint f(x,y) \Lambda(\Phi(x,y)) \,dx\,dy \qquad\text{(by the dominated convergence theorem)} \\
           = & \iint \Lambda(f(x,y) \Phi(x,y)) \,dx\,dy.
\end{align*}
Therefore, $\iint f(x,y) \Phi(x,y) \,dx\,dy = I$.
\end{pf}

\section{An interpolation theorem}

\begin{thm}\label{T:IT}
Let $1 \le p_0, p_1, p'_0, p'_1 \le \infty$ and suppose that
$T:L^{p_0}(\R^2) \cap L^{p_1}(\R^2) \to S^{p'_0} \cap S^{p'_1}$ is a linear
transformation which satisfies
\begin{align*}
\norm{Tf}_{S^{p'_0}} \le & M_0 \norm{f}_{p_0} \qquad \text{ and } \\
\norm{Tf}_{S^{p'_1}} \le & M_1 \norm{f}_{p_1}.
\end{align*}
Then for each $f \in L^{p_0}(\R^2) \cap L^{p_1}(\R^2)$ and each $t \in (0,1)$,
$Tf \in S^{p'_t}$ and 
$$ \norm{Tf}_{S^{p'_t}} \le M_t \norm{f}_{p_t},$$
where
\begin{align*}
\frac{1}{p_t} & = \frac{t}{p_1} + \frac{1-t}{p_0}, \\
\frac{1}{p'_t} & = \frac{t}{p'_1} + \frac{1-t}{p'_0} \qquad \text{ and } \\
M_t & = M_0^{1-t} M_1^t.
\end{align*}
\end{thm}

\begin{pf}
This follows immediately from the abstract (Calderon-Lions) interpolation
theorem once we know that $\{L^p(\R^2) | p \in [1,\infty]\}$ and 
$\{S^p | p \in [1,\infty]\}$ form complex interpolation scales.  For
this, see \cite[p38, Example 1 and p44, Proposition 8]{Reed-Simon-II}.
Note that we must take $S^\infty$ to be the space of compact operators with
the operator norm for all this to work.
\end{pf}

\section{Non-commutative H\"older inequality}

\begin{thm}\label{T:NCHI}
Let $1 \le p \le \infty$ and $p^{-1} + {p'}^{-1}=1$.  If $A\in S^p$ and
$B\in S^{p'}$, then $AB \in S^1$ and
$$\norm{AB}_{S^1} \le \norm{A}_{S^p} \norm{B}_{S^{p'}}.$$
\end{thm}

\begin{pf}
This follows from Theorem \ref{T:IT}.  For details,
see \cite{Reed-Simon-II}.
\end{pf}

\section{The inverse Alpha transform}

\begin{defn}
Let $f \in L^1(\R^2)$.  Then the inverse Alpha transform of $f$ is
the bounded operator defined by
$$ \Theta(f) = \iint f(x,y) M_{-y}T_{-x} \,dx \,dy $$
\end{defn}

Note that this integral exists by Lemma \ref{L:LOI1}.

\begin{lem}
If $f \in L^1(\R^2)$, then for all $g \in L^2(\R)$,
$$[\Theta(f)g](v) = \int K(v,w)g(w)\,dw \quad a.e.$$
where
$$K(v,w)=\int f(v-w, y) e^{-2\pi i yv} \,dy.$$
\end{lem}

\begin{pf}
For any $h \in L^2(\R)$,
\begin{align*}
\int [\Theta(f)g](v) \bar{h(v)} \,dv
      =& \int\(\iint f(x,y) M_{-y}T_{-x}\,dx\,dy \,g\)(v) \bar{h(v)} \,dv \\
      =& \iint f(x,y) \int (M_{-y} T_{-x}g)(v) \bar{h(v)} \,dv \,dx\,dy \\
       & \qquad \text{(by definition of integral)} \\
      =& \iint f(x,y) \int e^{-2\pi iyv} g(v-x) \bar{h(v)}\,dv \,dx\,dy \\
      =& \iiint f(x,y) e^{-2\pi iyv} \,dy\, g(v-x)\,dx\, \bar{h(v)}\,dv \\
       & \qquad \text{(by Fubini's theorem)} \\
      =& \iiint f(v-w,y) e^{-2\pi iyv}\,dy\, g(w)\,dw\, \bar{h(v)}\,dv \\
      =& \iint K(v,w) g(w) \,dw\, \bar{h(v)} \,dv.
\end{align*}
The application of Fubini's theorem is justified by the fact that
$F \in L^1(\R^3)$ if
$$ F(x,y,v)=f(x,y) e^{-2\pi iyv} g(v-x) \bar{h(v)}. \qed$$
\renewcommand{\qed}{}
\end{pf}

\section{A minimalist Alpha inversion formula}

\begin{lem}
If $X \in S^1$, then 
$$\norm{\alpha(X)}_\infty \le \norm{X}_{S^1}.$$
\end{lem}

\begin{pf}
For any $Y\in S^1$, we have $|\tr(Y)| \le \tr(|Y|)$
by the spectral theorem.  So for any $X\in S^1$,
\begin{align*}
\norm{\alpha(X)}_\infty
   =&   \sup_{(x,y)\in\R^2}|\tr(T_x M_y X)| \\
 \le&   \sup_{(x,y)\in\R^2}\tr(|T_x M_y X|) \\
   =&   \tr(|X|) \qquad\text{(since $T_x M_y$ is unitary)} \\
   =&   \norm{X}_{S^1}. \qed
\end{align*}
\renewcommand{\qed}{}
\end{pf}

Note that if $X$ is given by an integral kernel $K$ of Schwartz class, then
$$ \alpha(X)(x,y) = \int e^{2\pi i yv} K(v, v-x) \,dv. $$

\begin{thm}\label{T:SWIF1}
If $f$ is a Schwartz class function on $\R^2$, then
$$ \alpha(\Theta(f)) = f. $$
\end{thm}

\begin{pf}
By the Schwartz-Plancherel theorem for the classical Fourier transform,
the kernel $K$ of $\Theta(f)$ is of Schwartz class.  So
\begin{align*}
\alpha(\Theta(f))(x,y)
  =& \int e^{2\pi i yv} K(v, v-x) \,dv \\
  =& \int e^{2\pi i yv} \int f(x, y') e^{-2\pi i y'v} \,dy' \,dv \\
  =& f(x,y) \qquad\text{(by the classical Fourier inversion formula.)} \qed
\end{align*}
\renewcommand{\qed}{}
\end{pf}

\begin{cor}\label{C:SWIF2}
If $X \in S^1$ and $\alpha(X)$ is of Schwartz class, then
$$ X = \Theta(\alpha(X)). $$
\end{cor}

\begin{pf}
By Theorem \ref{T:SWIF1}, we have
$$ \alpha(\Theta(\alpha(X))) = \alpha(X). $$
But $\alpha$ is injective.  So
$$ \Theta(\alpha(X))=X. \qed$$
\renewcommand{\qed}{}
\end{pf}

\section{Non-commutative Riemann-Lebesgue Lemma}

\begin{thm}\label{T:NCRL}
If $f \in L^1(\R^2)$, then $\norm{\Theta(f)}_\infty \le \norm{f}_1$
and $\Theta(f) \in S^\infty$.
\end{thm}

\begin{pf}
Firstly,
\begin{align*}
\norm{\Theta(f)}_\infty
    =& \left\Vert \iint f(x,y) M_{-y} T_{-x} \,dx \,dy \right\Vert_\infty \\
 \le & \iint |f(x,y)| \norm{M_{-y} T_{-x}}_\infty \,dx \,dy \\
    =& \iint |f(x,y)| \,dx \,dy \\
    =& \norm{f}_1.
\end{align*}
Now, there is a sequence $\{f_k\}$ of compactly supported smooth functions such that $f_k \to f$
in $L^1(\R^2)$.  Moreover, for each $k$, the operator
$\Theta(f_k) \in S^{\infty}$.  But, by the previous calculation,
$\Theta(f_k) \to \Theta(f)$ in operator norm.  Since $S^{\infty}$
is closed in the operator norm, $\Theta(f) \in S^{\infty}$.
\end{pf}

\section{Non-commutative Plancherel theorem}

\begin{thm}\label{T:NCPT}
$\Theta$ extends to an isometry
$$\Theta:L^2(\R^2) \to S^2.$$
\end{thm}

\begin{pf}
Assume first that $f \in L^1(\R^2) \cap L^2(\R^2)$.  Then
$\Theta(f)$ is given by the kernel
$$ K(v,w) = \int f(v-w, y) e^{-2\pi iyv} \,dy.$$
So
\begin{align*}
\norm{\Theta(f)}_{S^2}^2
    =& \iint |K(v,w)|^2 \,dv \,dw \\
    =& \iint \left\vert \int f(v-w,y) e^{-2\pi iyv} \,dy 
       \right\vert^2 \,dv \,dw \\
    =& \iint |f(u,v)|^2 \,dv\,du\qquad \text{(by the classical Plancherel theorem)} \\
    =& \norm{f}_2.
\end{align*}
The rest is clear.
\end{pf}

\section{Non-commutative Hausdorff-Young theorem}

\begin{thm}\label{T:NCHY}
Let $1 \le p \le 2$ and $p^{-1} + {p'}^{-1} =1$.  Then $\Theta$ extends to
a bounded operator
$$\Theta:L^p(\R^2) \to S^{p'}.$$
\end{thm}

\begin{pf}
The endpoint estimates are given by Theorem \ref{T:NCRL}
and Theorem \ref{T:NCPT}.  The result now follows
from Theorem \ref{T:IT}.
\end{pf}

\section{An approximation lemma}

\begin{defn}
The {\em non-commutative Schwartz space} is the space $\SS$
of finite rank operators $X:L^2(\R) \to L^2(\R)$ such that
$$ X = \sum_{k=1}^n \phi_k \otimes \bar{\psi_k}\qquad \text{ with } 
\phi_k \in \S,$$
where $\S$ is the space of Schwartz class functions on $\R$.
\end{defn}

\begin{lem}\label{L:AL}
Let $X\in S^1$, $\tr(X)=0$ and $\e > 0$.  Then there exists
$Z \in \SS$ such that $\tr(Z)=0$ and 
$$\norm{X-Z}_{S^1} < \e.$$
\end{lem}

\begin{pf}
It is well known (see e.g. \cite{Reed-Simon}) that there exists a finite rank operator
$$X_1 = \sum_{k=1}^n \phi_k \otimes \bar{\psi_k}$$
such that
$\norm{\psi_k}_2 = 1$ and $\norm{X-X_1}_{S^1} < \frac{\e}{4}$.
It is also well known that there exist $\phi'_k \in \S$ such
that $\norm{\phi_k - \phi'_k}_2 < \frac{\e}{4n}$.  Set
$$X_2=\sum_{k=1}^n \phi'_k \otimes \bar{\psi_k}.$$
Then
\begin{align*}
\norm{X_1-X_2}  =& \left\Vert\sum_{k=1}^n (\phi_k-\phi'_k)\otimes\bar{\psi_k}\right\Vert_{S^1} \\
           \le& \sum_{k=1}^n \norm{\phi_k-\phi'_k}_2 \norm{\psi_k}_2 \\
             <& \sum_{k=1}^n \frac{\e}{4n} \\
             =& \frac{\e}{4}.
\end{align*}
Fix $W\in \SS$ such that $\norm{W}=\tr(W)=1$ and define $Z=X_2-\tr(X_2) W$. 
Then $Z\in \SS$, $\tr(Z)=0$ and
\begin{align*}
\norm{X_2-Z}_{S^1} \le& |\tr(X_2)| \\
                  =& \left|\sum_{k=1}^n \<\phi'_k,\psi_k\>\right| \\
                  =& \left|\sum_{k=1}^n \(\<\phi_k,\psi_k\>-\<\phi_k-
                                                   \phi'_k,\psi_k\>\)\right| \\
                \le& \left|\sum_{k=1}^n\<\phi_k,\psi_k\>\right|+
                             \sum_{k=1}^n |\<\phi_k-\phi'_k,\psi_k\>| \\
                \le& |\tr(X_1)| + 
                         \sum_{k=1}^n \norm{\phi_k-\phi'_k}_2\norm{\psi_k}_2\\
                  <& |\tr(X_1)| + \frac{\e}{4}.
\end{align*}
But 
\begin{align*}
|\tr(X_1)|   =& |\tr(X) - \tr(X-X_1)| \\
           \le& |\tr(X)| + |\tr(X-X_1)| \\
             =& |\tr(X-X_1)| \\
           \le& \tr(|X-X_1|) \\
             =& \norm{X-X_1}_{S^1} \\
             <& \frac{\e}{4}.
\end{align*}
Therefore, $\norm{X_2-Z}_{S^1} < \frac{\e}{2}$.
Therefore, $\norm{X-Z}<\e$.  
\end{pf}

\section{The action of differential operators}

\begin{defn}
For $\phi \in \S$, define
\begin{align*}
(P\phi)(t) = & \frac{d\phi}{dt} \\
(Q\phi)(t) = & 2 \pi i t \phi(t).
\end{align*}
\end{defn}

\begin{lem}
If $X \in \SS$, then $PX, QX \in \SS$.  Moreover, 
\begin{align*}
\alpha(PX) =& \(\frac{\partial}{\partial x} - 2\pi iy \)\alpha(X) \\
\alpha(QX) =& \frac{\partial}{\partial y} \alpha(X).
\end{align*}
\end{lem}

\begin{pf}
It follows immediately from the definition that $PX, QX \in \SS$.
Since for any $\phi \in \S$, 
$\lim_{h \to 0} \frac{T_h -I}{h}\phi = P\phi$ in the $L^2(\R)$ sense,
$\lim_{h \to 0} \frac{T_h -I}{h}X = PX$ in $S^1$-norm.  So
\begin{align*}
\frac{\partial}{\partial x} \alpha(X)(x,y)
    =& \frac{\partial}{\partial x} \tr(T_x M_y X) \\
    =& \lim_{h \to 0} \frac{\tr(T_{x+h} M_y X) - \tr(T_x M_y X)}{h} \\
    =& \lim_{h \to 0} 
           \tr\(\frac{e^{2\pi i yh} T_x M_y T_h X - T_x M_y X}{h}\) \\
    =& \lim_{h \to 0} \tr\(\frac{e^{2\pi i yh} T_x M_y T_h X -
            e^{2\pi i yh} T_x M_y X}{h} \right. \\
     & \qquad + \left. \frac{e^{2\pi i yh} T_x M_y X - T_x M_y X}{h} \) \\
    =& \lim_{h \to 0} \tr\(T_x M_y \frac{T_h -I}{h} X \) +
         \frac{e^{2\pi i yh} -1}{h} \tr(T_x M_y X) \\
    =& \tr(T_x M_y PX) + 2\pi i y \tr(T_x M_y X) \\
    =& (\alpha(PX) + 2\pi i y \alpha(X))(x,y).
\end{align*}
So
$$\alpha(PX) = \(\frac{\partial}{\partial x} - 2\pi i y\)\alpha(X). $$
By essentially the same argument, we get
$$\alpha(QX) = \frac{\partial}{\partial y} \alpha(X). \qed$$
\renewcommand{\qed}{}
\end{pf}

\begin{lem}
If $X \in \SS$ and $q \in \RA(\R^2)$, then $\im(q\cdot X) \subseteq \S$.
\end{lem}

\begin{pf}
Without loss of generality, we may assume $X=\phi\otimes\bar{\psi}$ with
$\phi \in \S$.  Let $g \in L^2(\R)$.  Since the map
$ (x,y) \mapsto \<g, T_{y_1} M_{-{x_1}}\psi\>$ is bounded and continuous,
$q({x_1},{y_1}) \<g, T_{y_1} M_{-{x_1}}\psi\> \in \RA(\R^2)$.
Since
$({x_1},{y_1}) \mapsto T_{y_1} M_{-{x_1}}\phi$ is a polynomially bounded continuous map $\R^2 \to \S$,
$$I = \iint q({x_1},{y_1}) \<g, T_{y_1} M_{-{x_1}} \psi\> T_{y_1}M_{-{x_1}}\phi \,d{x_1} \,d{y_1}$$
exists in $\S$ by Lemma \ref{L:LOI2}.  Now, for any $h \in L^2(\R)$,
\begin{align*}
\int I(v) \bar{h(v)} \,dv
    =& \iint q({x_1},{y_1}) \<(({x_1},{y_1})\cdot X)g, h\> \,d{x_1} \,d{y_1} 
           \qquad\text{(by definition of integral)} \\
    =& \<(q\cdot X)g,h\>.
\end{align*}
It follows that
$$ I = (q\cdot X)g \quad a.e. $$
In particular $(q\cdot X)g \in \S$.
\end{pf}

\begin{lem}\label{L:PqX}
Let $X \in \SS$ and $q \in \RA(\R^2)$.  Then
$$ P(q \cdot X) = (-2\pi ix_1 q)\cdot X + q \cdot (PX).$$
In particular, $P(q \cdot X) \in S^1$.
\end{lem}

\begin{pf}
Since for any $\phi \in \S$,
$$ \lim_{h \to 0} \(\frac{T_h - I}{h}\)\phi = P\phi$$
in the $L^2(\R)$-sense and $\im(q \cdot X) \subseteq \S$, we have
$$ P(q \cdot X) = \lim_{h \to 0} \frac{T_h - I}{h}(q \cdot X) $$
in the strong operator topology.
For the same reason,
$$ PX = \lim_{h \to 0} \frac{T_h -I}{h} X$$
in $S^1$-norm.
So
\begin{align*}
\lim_{h \to 0} \frac{T_h - I}{h} (q \cdot X)
    =& \lim_{h \to 0} \frac{T_h - I}{h} \iint
           q(x_1,y_1) (x_1,y_1) \cdot X \,dx_1 \,dy_1 \\
    =& \lim_{h \to 0} \iint q(x_1,y_1) \frac{T_h -I}{h}
           (x_1,y_1)\cdot X \,dx_1 \,dy_1 \\
    =& \lim_{h \to 0} \iint q(x_1,y_1) \frac1h 
           \(e^{-2\pi i x_1 h} T_{y_1} M_{-x_1} T_h X M_{x_1} T_{-y_1} \right.\\
     & \qquad\left. - T_{y_1} M_{-x_1} X M_{x_1} T_{-y_1}\) \,dx_1 \,dy_1 \\
    =& \lim_{h \to 0} \iint q(x_1,y_1)
       \(\frac{e^{-2\pi i x_1 h} -1}{h} T_{y_1} M_{-x_1} T_h X M_{x_1} T_{-y_1}\right.\\
     & \qquad\left. -T_{y_1}M_{x_1} \frac{T_h -I}{h}X M_{x_1} T_{-y_1}\) \,dx_1 \,dy_1 \\
    =& \iint q(x_1,y_1)(-2\pi i x_1 (x_1,y_1)\cdot X + (x_1,y_1)\cdot (PX))\,dx_1\,dy_1\\
    =& (-2\pi i x_1 q)\cdot X + q \cdot (PX)
\end{align*}
in $S^1$-norm.
This proves the claim.
\end{pf}

\begin{lem}\label{L:QqX}
Let $X \in \SS$ and $q \in \RA(\R^2)$.  Then
$$ Q(q \cdot X) = (-2\pi iy_1 q)\cdot X + q \cdot (QX).$$
In particular, $Q(q \cdot X) \in S^1$.
\end{lem}

\begin{pf}
This is proved by the same sort of reasoning as Lemma \ref{L:PqX}.
\end{pf}

\begin{lem}\label{L:alphaPqX}
If $X \in \SS$ and $q\in \RA(\R^2)$, then
$$ \alpha(P(q \cdot X))=
   \(\frac{\partial}{\partial x} - 2\pi iy \) \alpha(q \cdot X). $$
\end{lem}

\begin{pf}
By Lemma \ref{L:PqX}, we have
$$ P(q\cdot X) = (-2\pi ix_1 q)\cdot X + q \cdot (PX).$$
It follows that
\begin{align*}
\alpha(P(q \cdot X))
    =& \frac{\partial\hat q}{\partial x} \alpha(X) +
            \hat q \alpha(PX) \\
    =& \frac{\partial\hat q}{\partial x} \alpha(X) +
            \hat q \(\frac{\partial}{\partial x} - 2\pi i y\) \alpha(X) \\
    =& \frac{\partial}{\partial x}(\hat q \alpha(X)) - 2\pi i y \alpha(X) \\
    =& \(\frac{\partial}{\partial x} - 2\pi iy \)(\hat q \alpha(X)) \\
    =& \(\frac{\partial}{\partial x} - 2\pi iy\) \alpha(q \cdot X). \qed
\end{align*}
\renewcommand{\qed}{}
\end{pf}

\begin{lem}\label{L:alphaQqX}
If $X \in \SS$ and $q\in \RA(\R^2)$, then
$$ \alpha(Q(q \cdot X))=
     \frac{\partial}{\partial y}\alpha(q \cdot X). $$
\end{lem}

\begin{pf}
This is proved in the same way as Lemma \ref{L:alphaPqX}.
\end{pf}

\section{The harmonic oscillator}\label{S:HO}

\begin{defn}
The {\em harmonic oscillator} is the differential operator
$$H= P^2+Q^2.$$
\end{defn}

\begin{lem}\label{L:HO}
If $X \in \SS$ and $q \in \RA(\R^2)$, then $H(q\cdot X) \in S^1$ and
$$\alpha(H(q\cdot X)) = {\mathcal D}\alpha(q \cdot X),$$
where
$$ {\mathcal D} =\(\frac{\partial}{\partial x} -2\pi iy \)^2 +
                \(\frac{\partial}{\partial y} \)^2. $$
\end{lem}

\begin{pf}
By Lemma \ref{L:PqX} and Lemma \ref{L:QqX}, we have
\begin{align*}
P^2(q\cdot X) =& P((-2\pi i x_1 q)\cdot X) +P(q\cdot (PX))\qquad\text{and}\\
Q^2(q\cdot X) =& Q((-2\pi i y_1 q)\cdot X) +Q(q\cdot (QX)).
\end{align*}
So $P^2(q\cdot X)$, $Q^2(q\cdot X) \in S^1$ and hence $H(q\cdot X) \in S^1$.

By Lemma \ref{L:alphaPqX} and Lemma \ref{L:alphaQqX}, we have
\begin{align*}
\alpha(P^2(q\cdot X))
=& \alpha(P((-2\pi i x_1 q)\cdot X)) +\alpha(P(q\cdot(PX))) \\
=& \(\frac{\partial}{\partial x} - 2\pi iy\) \alpha((-2\pi ix_1 q)\cdot X)
   +\(\frac{\partial}{\partial x} - 2\pi iy\) \alpha(q\cdot (PX)) \\
=& \(\frac{\partial}{\partial x} - 2\pi iy\)
     (\alpha((-2\pi ix_1 q)\cdot X)+\alpha(q\cdot (PX))) \\
=& \(\frac{\partial}{\partial x} - 2\pi iy\)
     \(\frac{\partial \hat q}{\partial x} \alpha(X) +
      \hat q \(\frac{\partial}{\partial x} - 2\pi iy\) \alpha(X)\) \\
=& \(\frac{\partial}{\partial x} - 2\pi iy\)
     \(\frac{\partial}{\partial x}(\hat q \alpha(X)) - 2\pi iy \alpha(X)\) \\
=& \(\frac{\partial}{\partial x} - 2\pi iy\)^2(\hat q \alpha(X)) \\
=& \(\frac{\partial}{\partial x} - 2\pi iy\)^2 \alpha(q\cdot X).
\end{align*}
and
\begin{align*}
\alpha(Q^2(q\cdot X))
=& \frac{\partial}{\partial y}\alpha((-2\pi i y_1 q)\cdot X) +
   \frac{\partial}{\partial y}\alpha(q\cdot (QX)) \\
=& \frac{\partial}{\partial y}
    (\alpha((-2\pi i y_1 q)\cdot X) + \alpha(q\cdot (QX))) \\
=& \frac{\partial}{\partial y}
    \(\frac{\partial\hat q}{\partial y} \alpha(X) +
      \hat q \frac{\partial}{\partial y} \alpha(X)\) \\
=& \frac{\partial^2}{\partial y^2}(\hat q \alpha(X)) \\
=& \frac{\partial^2}{\partial y^2}\alpha(q\cdot X).
\end{align*}
The result now follows.
\end{pf}

\begin{thm}
There is a complete orthonormal set $\{\phi_k\}$ in $L^2(\R)$ such that
$$ H\phi_k = -2\pi(2k+1) \phi_k. $$
\end{thm}

\begin{pf}
See \cite[Lemma 10.34]{Roe}.
\end{pf}

\begin{cor}\label{C:p-series}
For $p > 1$,
$$ H^{-1} \in S^p.$$
\end{cor}

\begin{pf}
This follows from the $p$-series test.
\end{pf}

\section{A versal constant}

For $\delta > 0$, let $B_\delta=\{(x,y) \in \R^2 : \sqrt{x^2+y^2} < \delta\}$.  Fix
a radial smooth function
$\tau$ with support in $B_1$ and identically $1$
in a neighborhood of $0$.  Set
$$ V = \norm{\check\tau}_1.$$

Set
$$ \tau_\delta(x,y) = \tau(x/\delta, y/\delta).$$
Then
\begin{align*}
\check\tau_\delta(x_1, y_1)
    &= \iint e^{2\pi i (xx_1+yy_1)} \tau_\delta(x,y) \,dx\,dy \\
    &= \iint e^{2\pi i (xx_1+yy_1)} \tau(x/\delta, y/\delta) \,dx\,dy \\
    &= \iint e^{2\pi i \delta (xx_1+yy_1)} \tau(x,y) \delta^2 \,dx\,dy \\
    &= \delta^2 \check\tau(\delta x_1, \delta y_1).
\end{align*}
So for any $\delta > 0$,
\begin{align*}
\norm{\check\tau_\delta}_1
    &= \delta^2 \iint |\check\tau(\delta x_1, \delta y_1)| \,dx_1\,dy_1 \\
    &= \delta^2 \iint |\check\tau(x_1, y_1)| \delta^{-2} \,dx_1\,dy_1 \\
    &= \norm{\check\tau}_1 \\
    &= V.
\end{align*}

\section{The heart of the matter}

\begin{lem}\label{L:HOM}
If $X \in \SS$, $\tr(X)=0$ and $\e > 0$, then there exists $\p \in L^1(\R^2)$
with $\hat\p=1$ on a neighborhood of $(0,0)$ such that
$$\norm{\p \cdot X}_{S_1} < \e.$$
\end{lem}

\begin{pf}
Let $r=\sqrt{x^2+y^2}$. Let
$$ C_1  = \sup_{B_1} |\nabla\alpha(X)|. $$
Since $\alpha(X)$ is smooth and $\alpha(X)(0,0)=0$, it follows
from the mean value theorem that
$$ |\alpha(X)(x,y)| \le C_1 r \qquad\text{on $B_1$}.$$
Let
\begin{align*}
C_2  =& \sup_{B_1} |\Delta\alpha(X)| \\
D_1  =& \sup_{B_1} |\tau| \\
D_2  =& \sup_{B_1} |\nabla\tau| \\
D_3  =& \sup_{B_1} |\Delta\tau|.
\end{align*}
Then
\begin{align*}
{\mathcal D}(\tau_\delta \alpha(X)) 
    =& \(\Delta -4\pi iy \frac{\partial}{\partial x}- 4\pi^2y^2\) 
         (\tau_\delta \alpha(X)) \\
    =& (\Delta\tau_\delta)\alpha(X)+2\nabla\tau_\delta\cdot\nabla\alpha(X)+
         \tau_\delta \Delta\alpha(X) \\
     & -4\pi iy \(\frac{\partial\tau_\delta}{\partial x} \alpha(X)+
           \tau_\delta \frac{\partial}{\partial x} \alpha(X)\)\\
     & -4\pi^2y^2 \tau_\delta \alpha(X).
\end{align*}
Therefore, on $B_\delta$, we have by the Cauchy--Schwarz inequality that
\begin{align*}
|{\mathcal D}(\tau_\delta \alpha(X))|
  \le & [|\Delta\tau_\delta||\alpha(X)|+2|\nabla\tau_\delta||\nabla\alpha(X)|+
            |\tau_\delta| |\Delta\alpha(X)|  \\
      & +4\pi |y| \(\left|\frac{\partial\tau_\delta}{\partial x}\right||\alpha(X)|+
            |\tau_\delta| \left|\frac{\partial}{\partial x} \alpha(X)\right|\) \\
      & +4\pi^2|y|^2 |\tau_\delta| |\alpha(X)|]  \\
  \le & (|\Delta\tau_\delta||\alpha(X)| + 2|\nabla\tau_\delta||\nabla\alpha(X)|+
            |\tau_\delta| |\Delta\alpha(X)| \\
      & +4\pi r(|\nabla\tau_\delta||\alpha(X)|+ |\tau_\delta||\nabla\alpha(X)|)+
            4\pi^2r^2|\tau_\delta||\alpha(X)|)  \\
  \le & (\delta^{-2} D_3 C_1 r + 2 \delta^{-1} D_2 C_1 + D_1 C_2 \\
      & +4\pi r (\delta^{-1} D_2 C_1 r + D_1 C_1) +
         4\pi^2 r^2 D_1 C_1 r)  \\
  \le & (\delta^{-2} D_3 C_1 r + \delta^{-1}(2 D_2 C_1 + 4\pi D_2 C_1)+
        (D_1 C_2 + 4\pi D_1 C_1 + 4\pi^2 D_1 C_1))  \\
  \le & (A_1 \delta^{-2} r + A_2 \delta^{-1} + A_3).
\end{align*}
Moreover, ${\mathcal D}(\tau_\delta \alpha(X))$ is supported in $B_\delta$.
Fix $p \in (1,2)$.  Then
\begin{align*}
\norm{{\mathcal D}(\tau_\delta \alpha(X))}_p
  \le & A_1 \delta^{-2} \(\iint_{B_\delta} r^p \,dx\,dy \)^{1/p}+
        A_2 \delta^{-1} \(\iint_{B_\delta} \,dx\,dy \)^{1/p} \\
      & + A_3 \(\iint_{B_\delta} \,dx\,dy \)^{1/p} \\
  \le & (2\pi)^{1/p}
        \[A_1 \delta^{-2} \(\frac{\delta^{p+2}}{p+2}\)^{1/p}+
          A_2 \delta^{-1} \(\frac{\delta^2}{2}\)^{1/p} +
          A_3 \(\frac{\delta^2}{2}\)^{1/p}\] \\
     =& (2\pi)^{1/p}
        \[(A_1 (p+2)^{-1/p} + A_2 2^{-1/p}\]\delta^{2/p-1} + A_3 2^{-1/p} \delta^{2/p}.
\end{align*}
Therefore,
$$\lim_{\delta \to 0} {\mathcal D}(\tau_\delta \alpha(X)) = 0 \qquad\text{in $L^p(\R^2)$}.$$

If $p^{-1}+{p'}^{-1}=1$, then we have
\begin{align*}
\lim_{\delta \to 0} H(\check\tau_\delta \cdot X) 
=& \lim_{\delta \to 0}\Theta \[{\mathcal D} (\tau_\delta \alpha(X))\]
\qquad\text{(by Corollary \ref{C:SWIF2} and Lemma \ref{L:HO})} \\
=& \Theta\[\lim_{\delta \to 0} {\mathcal D} (\tau_\delta \alpha(X))\]
\qquad\text{(by Theorem \ref{T:NCHY})} \\
=& 0 \qquad\text{in $S^{p'}$}.
\end{align*}
By Corollary \ref{C:p-series}, $H^{-1} \in S^p$.  So by Theorem \ref{T:NCHI}, we have
\begin{align*}
\lim_{\delta \to 0}\check \tau_\delta \cdot X 
      =& H^{-1}(\lim_{\delta \to 0} H(\check \tau_\delta \cdot X)) \\
      =& 0 \qquad\text{in $S^1$}.
\end{align*}
Therefore, there exists $\delta_0>0$ such that 
$\norm{\check\tau_{\delta_0} \cdot X}_{S^1} < \e$.
Take $\p=\check \tau_{\delta_0}$.  Then $\hat\p=1$ in a neighborhood of $(0,0)$
and $\norm{\p\cdot X}_{S^1} < \e$.
\end{pf}

\begin{PMT*}
Now, if $X \in S^1$ and $\tr(X)=0$, by Theorem \ref{L:AL} we can find $X' \in \SS$ such that
$\norm{X-X'}_{S^1} < \frac{\e}{2V}$ and $\tr(X')=0$.  Then by Lemma \ref{L:HOM} we can find
$\p \in L^1(\R^2)$ such that $\norm{\p \cdot X'}_{S^1}<\frac{\e}{2}$.  Thus
\begin{align*}
\norm{\p \cdot X}_{S^1} =& \norm{\p \cdot [X' - (X'-X)]}_{S^1} \\
                  \le & \norm{\p \cdot X'}_{S^1} +
                        \norm{\p \cdot (X'-X)}_{S^1} \\
                  \le & \frac{\e}{2} + \norm{\p}_1 \norm{X-X'}_{S^1} \\
                    < & \e. \qed
\end{align*}
\end{PMT*}

\bibliographystyle{amsplain}
\bibliography{nsi-v4}

\end{document}